\documentstyle[12pt]{article}
\def\R{\hbox{{\rm I}\kern-0.2em{\rm R}\kern0.2em}}
\def\D{\hbox{{\rm I}\kern-0.2em{\rm D}\kern0.2em}}

\def\be{\begin{equation}}
\def\ee{\end{equation}}

\def\({\left(}
\def\){\right)}
\def\[{\left[}
\def\]{\right]}
\def\bc{\begin{center}}
\def\ec{\end{center}}

\begin{document}

{\bf\Large \textbf{Reduction of fourth order ordinary differential\\
equations to second and third order Lie\\
linearizable forms}}
\begin{center}

\textbf{Hina M. Dutt, Asghar Qadir} \end{center}
School of Natural Sciences\\
National University of Sciences and Technology\\
Campus H-12, 44000, Islamabad, Pakistan\\
E-mail: hinadutt@yahoo.com; aqadirmath@yahoo.com

{\bf Abstract}.\\
Meleshko presented a new method for reducing third order autonomous
ordinary differential equations (ODEs) to Lie linearizable second
order ODEs. We extended his work by reducing fourth order autonomous
ODEs to second and third order linearizable ODEs and then applying
the Ibragimov and Meleshko linearization test for the obtained ODEs.
The application of the algorithm to several ODEs is also presented.

\section{Introduction}

First order ODEs can always be linearized \cite{fm1} by point
transformations \cite{l1}. Lie \cite{l2} showed that all
linearizable second order ODEs must be cubically semi-linear, i.e.
\begin{eqnarray}
y''+ a_1(x,y)y'^3 - a_2(x,y)y'^2 + a_3(x,y)y' - a_4(x,y)=0~,\label{secondlform}
\end{eqnarray}
the coefficients $a_1, a_2, a_3, a_4$ satisfy an over-determined
integrable system of four constraints involving two auxiliary
functions, which Tresse wrote in more usable form \cite{t}
\begin{eqnarray}
3(a_1a_3)_x - 3a_4a_{1y} - 6a_1a_{4y} - 2a_2a_{2x} + a_2a_{3y} -
3a_{1xx} +2a_{2xy} -
a_{3yy}=0~,\nonumber\\
3(a_4a_2)_y - 3a_1a_{4x} - 6a_4a_{1x} - 2a_3a_{3y} + a_3a_{2x}
+3a_{4yy} - 2a_{3xy} + a_{2xx}=0~.\label{secondlc}
\end{eqnarray}
We call such equations \emph{Lie linearizable}.\\
\\ Chern \cite{c1, c2} and Grebot \cite{g1, g2} extended the linearization programme to
the third order using contact and point transformations respectively
to obtain linearizability criteria for equations reducible to the
forms $u'''(t)=0$ and $u'''(t)+u(t)=0.$ It was shown \cite{mL} that
there are three classes of third order ODEs that are linearizable by
point transformations, viz. those that reduce to the above two forms
or $u'''(t)+\alpha(t)u(t)=0.$ Neut and Petitot \cite{np} dealt with the general third order ODEs.
Ibragimov and Meleshko (IM) \cite{im} used the original Lie
procedure \cite{l2} of point transformation to determine the linearizability criteria for third order ODEs. They showed that any third order ODE $y'''=f(x,y,y',y'')$ obtained from a linear equation $u''' + \alpha(t)u=0$ by means of
point transformations $t=\varphi(x,y),u=\psi(x,y),$
must belong to one of the following two types of equations.\\
\textbf{Type I}: If $\varphi_y=0$ the equations that are
linearizable are of the form
\begin{eqnarray}
y''' +(a_1y' + a_0)y'' + b_3y'^3 + b_2y'^2 + b_1y' + b_0=0~.\label{thirdtype1}
\end{eqnarray}
\textbf{Type II}: If $\varphi_y\neq0$ , set
$r(x,y)=\varphi_x/\varphi_y$~, equations are of the form
\begin{eqnarray}
y'''+{1\over y'+r}[-3(y'')^2 + (c_2y'^2+c_1y'+c_0)y''\nonumber\\
+d_5y'^5+d_4y'^4+d_3y'^3+d_2y'^2+d_1y'+d_0]=0~,\label{thirdtype2}
\end{eqnarray}
where all coefficients $a_i,b_i,c_i,d_i$, being the functions of $x$ and $y$, satisfy certain constraint requirements. Afterwards Ibragimov, Meleshko and Suksern \cite{ims1, ims2} used the point and contact transformations to
determine the criteria for the linearizability of fourth order
scalar ODEs. Meleshko \cite{M} provided a simple algorithm to reduce
third order ODEs of the form $y'''=f(y,y',y'')$ to second order
ODEs. If the reduced equations satisfy the Lie
linearizability criteria, they can then be solved by linearization. Meleshko showed that a third order ODE is reducible to the second order linearizable ODE if it is of the form
\begin{eqnarray}
y'''+A(y,y')y''^3+B(y,y')y''^2+c(y,y')y''+D(y,y')~,
\end{eqnarray}
where the coefficients $A,B,C,D$ satisfy certain constraints.
\\\\
In the present paper we extend Meleshko's procedure to the fourth
order ODEs in the cases that the equations do not depend explicitly
on the independent or the dependent variable (or both) to reduce it
to third (respectively second) order equations. Once the order is
reduced we can apply the IM (or Lie) linearization test. If the
reduced third (or second) order ODE satisfies the IM (or Lie)
linearization test, then after finding a linearizing transformation,
the general solution of the original equation is obtained by
quadrature. So this method is effective in the sense that it reduces
many ODEs, that cannot be linearized, to lower order linearizable
forms. This is one of the motivations for studying this method.
Another hope for the study of the linearization problem is that by
using it we may be able to provide a complete classification of ODEs
according to the number of arbitrary initial conditions that
can be satisfied \cite{mq2}.\\
\\

\section{Equations reducible to linearizable forms}
Meleshko had only treated the special case of independence of $x$
for third order ODE. We include independence of $y$ for completeness
before proceeding to the fourth order.\\\\
\textbf{Third order ODEs independent of $y$}\\\\
Taking $y'$ as the independent variable $u(x)$, we convert the ODE
\begin{eqnarray} y'''=f(x,y',y'')~,
\end{eqnarray} to the second order ODE
\begin{eqnarray}
u''=f(x,u,u')~,\label{thirdindy-general}
\end{eqnarray}
which is linearizable by Lie's criteria if it is cubically
semi-linear with the coefficients satisfying conditions (\ref{secondlc}).\\
\\
Hence (\ref{thirdindy-general}) is reducible to second order linearizable form if and
only if
\begin{eqnarray}
f(x,y',y'')=-c(x,y')y''^3+g(x,y')y''^2-h(x,y')y''+d(x,y')~,\label{thirdindy-lform}
\end{eqnarray}
with the coefficients satisfying
\begin{eqnarray}
3(ch)_x-3dc_{y'}-6cd_{y'}-2gg_x+gh_{y'}-3c_{xx}+2g_{xy'}-h_{y'y'}=0~,\nonumber\\
3(dg)_{y'}-3cd_x-6dc_x-2hh_{y'}+hg_x+3d_{y'y'}-2h_{xy'}+g_{xx}=0~.\label{thirdindy-lc}
\end{eqnarray}
\\
\textbf{Fourth order ODEs independent of $y$}\\\\
Since the variable $y$ is missing, by taking $y'$ as the new
dependent variable $u(x)$, the ODE
\begin{eqnarray}
y^{(4)}=f(x,y',y'',y''')~,\label{fourthindy-general}
\end{eqnarray}
is reduced to third order ODE \begin{eqnarray} u'''=f(x,u,u',u'')~.\label{reduced10}
\end{eqnarray}
Equation (\ref{reduced10}) is linearizable for the type I of Ibragimov and
Meleshko's criteria if and only if
\begin{eqnarray}
f(x,y',y'',y''')=-(a_1y''+a_0)y'''-b_3y''^3-b_2y''^2-b_1y''-b_0~,
\end{eqnarray} with the coefficients $a_i=a_i(x,y')~,~(i=0,1)$ and $b_j=b_j(x,y')~,~(j=0,1,2,3)$~,\label{fourthindy-lform}\\
\\
satisfying the conditions
\begin{eqnarray}
a_{0y'}-a_{1x}=0~,\nonumber\\
(3b_1-{a_0}^2-3a_{0x})_{y'}=0~,\nonumber\\
3a_{1x}+a_0a_1-3b_2=0~,\nonumber\\
3a_{1y'}+{a_1}^2-9b_3=0~,\nonumber\\
(9b_1-6a_{0x}-2{a_0}^2)a_{1x}+9(b_{1x}-a_1b_0)_{y'}+3b_{1y'}a_0-27b_{0y'y'}=0~.\label{fourthindy-lc}
\end{eqnarray}\\
Also the necessary and sufficient conditions for (\ref{reduced10}) to be
linearizable for the type II of Ibragimov and Meleshko's criteria
are
\begin{eqnarray}
f(x,y',y'',y''')={-1\over y''+r}[-3(y''')^2 + (c_2y''^2+c_1y''+c_0)y'''\nonumber\\
+d_5y''^5+d_4y''^4+d_3y''^3+d_2y''^2+d_1y''+d_0]~,\label{fourthindy-lform-type2}
\end{eqnarray}
and the coefficients $c_i=c_i(x,y')~,~(i=0,1,2)$ ,
$d_j=d_j(x,y')~,~(j=0,1,2,3,4,5)$ and $r=r(x,y')$ have to satisfy constraint equations which can be produced simply by replacing $y$ by $y'$ for the type II constraint equations in \cite{im}.\\\\
\textbf{Fourth order ODEs independent of $x$}\\\\
 The transformation $y'=u(y)$ will transform autonomous ODE of the
fourth order
\begin{eqnarray}
y^{(4)}=f(y,y',y'',y''')~,\label{fourthindx-general}
\end{eqnarray}
into the equation
\begin{eqnarray}
u^3u'''+4u^2u'u''+uu'^3-f(y,u,uu',u^2u''+uu'^2)=0~,\label{reduced1.1}
\end{eqnarray}
which is a third order ODE in $(y,u)$. It is linearizable by
Ibragimov Meleshko's criteria if it is of the form (\ref{thirdtype1}) i.e,
\begin{eqnarray}
f(y,u,u'u,u''u^2+uu'^2)=-u^3[(a_1u'+a_0)u''+b_3u'^3+b_2u'^2+b_1u'+b_0]\nonumber\\
+4u^2u'u''+uu'^3~,\label{reduced2}
\end{eqnarray}
where $a_i=a_i(y,u)~,~(i=0,1)$ and $b_j=b_j(y,u)~,~(j=0,1,2,3)~.$
With this (\ref{reduced1.1}) takes the form
\begin{eqnarray}
u'''+(a_1u'+a_0)u''+b_3u'^3+b_2u'^2+b_1u'+b_0=0~.\label{reduced3}
\end{eqnarray}
Transforming (\ref{reduced3}) into a fourth order ODE with $x$ as independent
variable and $y$ as dependent variable:
\begin{eqnarray}
y^{(4)}+(A_1y''+A_0)y'''+B_3y''^3+B_2y''^2+B_1y''+B_0=0~,\label{fourthindx-lform-type1}
\end{eqnarray}
where \begin{eqnarray} A_i=A_i(y,y')~,~ (i=0,1)~;~\quad
B_j=B_j(y,y')~, (j=0,1,2,3) \end{eqnarray} subject to the
identification of coefficients
\begin{eqnarray}
a_1=A_1+{4\over y'}~,~\quad a_0={A_0\over y'}~,~\quad
b_3=B_3+{A_1\over
y'}+{1\over y'^2}~, \nonumber\\
b_2={B_2\over y'}+{A_0\over y'^2}~,~\quad b_1={B_1\over
y'^2}~,~\quad b_0={B_0\over y'^3}~,
\end{eqnarray}
with the constraints
\begin{eqnarray}
y'^2A_{1y}-y'A_{0y'}+A_0=0~,\nonumber\\
y'^2(-3A_{0yy'})+y'(3B_{1y'}+3A_{0y}-2A_0A_{0y'})+(-6B_1+2A_0^2)=0~,\nonumber\\
y'^2(3A_{1y})+y'(A_0A_1-3B_2)+A_0=0~,\nonumber\\
y'^2(3A_{1y'}-9B_3+A_1^2)-y'A_1-5=0~,\nonumber\\
y'^4(-6A_{0y}A_{1y})+y'^3(9B_1A_{1y}-2A_0^2A_{1y}+9B_{1yy'})+y'^2(-18B_{1y}-9A_1B_{0y'}\nonumber\\
-9B_0A_{1y'}
+3A_0B_{1y'}-27B_{0y'y'})+y'(27A_1B_0-6A_0B_1+126B_{0y'})\nonumber\\
-180B_0=0~.\label{fourthindx-lc-type2}
\end{eqnarray}
\\
Also in order to make (\ref{reduced1.1}) linearizable of type II of Ibragimov and
Meleshko's criteria we have to take
\begin{eqnarray} f({y,u,uu',u^2u''+uu'^2})=-{u^3\over u'+r}[-3(u'')^2+(c_2u'^2+c_1u'+c_0)u''\nonumber\\
+d_5u'^5+d_4u'^4+d_3u'^3+d_2u'^2+d_1u'+d_0]+4u^2u'u''+uu'^3~,\label{reduced4}
\end{eqnarray}
where $c_i=c_i(y,u)~,~(i=0,1,2)$~,~ $d_j=d_j(y,u)~,~(j=0,1,2,3,4,5)$ and $r=r(y,u)~.$\\
\\
Considering the form (\ref{reduced4}) and converting (\ref{reduced1.1}) into fourth order with
$x$ as independent and $y$ as dependent variable, we have
\begin{eqnarray} y^{(4)} +{1\over y''+r_0}[-3(y''')^2 + (C_2y''^2+C_1y''+C_0)y'''\nonumber\\
+D_5y''^5+D_4y''^4+D_3y''^3+D_2y''^2+D_1y''+D_0]=0~,\label{fourthindx-lform-type2}
\end{eqnarray}
where \begin{eqnarray} C_i=C_i(y,y')~,~ (i=0,1,2)~;\quad
D_j=D_j(y,y')~,~ (j=0,1,2,3,4,5)~;~\quad r_0=r_0(y,y')~,\nonumber
\end{eqnarray} subject to the identification of coefficients
\begin{eqnarray}
c_2=C_2-{2\over y'}~,~\quad c_1=C_1+{4r_0\over y'}~,~\quad
c_0={C_0\over
y'^2}~,~\quad d_5={D_5\over y'^5}~, \nonumber\\
d_4=D_4+{C_2\over y'}-{2\over y'^2}~,~\quad d_3={D_3\over
y'}+{C_1\over y'}+{4r_0\over y'^2}-{3r_0\over y'^3}~,\nonumber\\
d_2={D_2\over y'^2}+{C_0\over y'^3}~,~\quad d_1={D_1\over
y'^3}~,~\quad d_0={D_0\over y'^4}~,~\quad r={r_0\over y'}~,
\end{eqnarray}
with the constraints (\ref{C1})--(\ref{C9}) (presented in the appendix).\\
\\
\textbf{Fourth order ODEs independent of $x$ and $y$}\\\\
By considering $y'$ as independent and $y''$ as dependent variable,
we convert the equation \begin{eqnarray} y^{(4)}=f(y',y'',y''')~,\label{fourthindxy-general}
\end{eqnarray} into a second order ODE:
\begin{eqnarray}
u^2u''+uu'^2=f(y',u,uu')~.\label{reduced5}
\end{eqnarray}
For (\ref{reduced5}) to be Lie-linearizable we must have
\begin{eqnarray}
f(y',u,uu')=-u^2[A(y',u)u'^3+B(y',u)u'^2+C(y',u)u'+D(y',u)]+uu'^2~.\label{reduced6}
\end{eqnarray}
Hence (\ref{fourthindxy-general}) takes the form \begin{eqnarray}
y^{(4)}+a(y',y'')y'''^3+b(y',y'')y'''^2+c(y',y'')y'''+d(y',y'')=0~,\label{fourthindxy-lform}
\end{eqnarray}
where $a,b,c$ and $d$ must satisfy the constraints:
\begin{eqnarray}
(3a_{y'y'})y''^4+(2bb_{y'}-3ca_{y'}-3ac_{y'}-2b_{y'y''})y''^3+(2b_{y'}-bc_{y''}\nonumber\\
+3a_{y''}d+6ad_{y''}-c_{y''y''})y''^2+(bc-9ad-3c_{y''})y''-c=0~,\nonumber\\
(by'y')y''^4+(b_{y'}c+3d_{y''}b-3d_{y'}a-6a_{y'}d-2c_{y'y''})y''^3+(c_{y'}+3d_{y''}\nonumber\\
-6bd+3b_{y''}d-2cc_{y''}+3d_{y''y''})y''^2+(2c^2-6d-12d_{y''})y''+15d=0~.\label{fourthindxy-lc}
\end{eqnarray}

Thus we have the following theorems.\\
\\ {\bf Theorem 1.} \emph{Equation} (\ref{fourthindx-lform-type1})
\emph{is reduced to the third order linearizable form if and only if
it obeys} (\ref{fourthindx-lc-type2}).
\\\\
{\bf Theorem 2.} \textit{Equation} (\ref{fourthindx-lform-type2}) \emph{is reduced to the
third order linearizable form if
and only if it obeys} (\ref{C1})--(\ref{C9}) (presented in the appendix II).\\
\\
\textbf{Theorem 3}. \emph{Equation} (\ref{fourthindxy-lform}) \emph{is reduced to the
second order linearizable form if and only if it obeys} (\ref{fourthindxy-lc}).\\\\
\textbf{Remark}: If we have a fourth order ODE of the form
\begin{eqnarray}
y^{(4)}=-f(x,y)y'^5+10{y''y'''\over y'}-15{y''^3\over y'^2}~,
\end{eqnarray}
with $f(x,y)$ linear in $x$, then we can convert it to a linear ODE
$x^{(4)}=f(x,y)$ by simply taking $x$ as dependent and $y$ as independent variables.

\section{Illustrative Examples}

{\bf Example 1.} The nonlinear fourth order ODE
\begin{eqnarray} y'y^{(4)}-y''y'''-3y'^2y'''+2y'^3y''+3y'^5=0~,\label{example1}
\end{eqnarray}
cannot be linearized by point or contact transformation. It has the
form (\ref{fourthindx-lform-type1}) with the coefficients $A_1=-1/y',$ $A_0=-3y',$ $B_3=B_2=0,$ $B_1=2y'^2,$ $B_0=3y'^5.$
One can verify that these coefficients satisfy the conditions
(\ref{fourthindx-lc-type2}). The transformation $y'=u(y)$ will reduce this ODE to the 3rd
order linearizable ODE
\begin{eqnarray}
u'''+{3\over u}u'u''-3u''-{3\over u}u'^2+2u'+3u=0~.\label{reduced7} \end{eqnarray}
By using transformation equations in \cite{im}, we arrive at the transformation
$t=e^y,$ $s=u^2$
which maps (\ref{reduced7}) to the linear third order
ODE $s'''-{2\over t^3}s=0$,
whose solution is given by $s=c_1t^{-1}+t^2\{c_2\cos(\sqrt{2}\ln t)+c_3\sin (\sqrt{2} \ln t)\}$,
where $c_i$ are arbitrary constants. By using the above transformation we get the solution of (\ref{reduced7}) given by
$u=\pm \sqrt{c_1e^{-y}+e^{2y}\{c_2\cos (\sqrt{2}y)+c_3\sin (\sqrt{2}y)\}}$. Hence the general solution of (\ref{example1}) is obtained by
taking the quadrature
\begin{eqnarray}
\int {dy\over \sqrt{c_1e^{-y}+e^{2y}\{c_2\cos (\sqrt{2}y)+c_3\sin (\sqrt{2}y)\}}}=\pm x+c_4~,
\end{eqnarray}
where $c_i$ are arbitrary constants.\\
\\
{\bf Example 2.} The nonlinear ODE
\begin{eqnarray}
y^2y'^2y^{(4)}-10y^2y'y''y'''-3yy'^3y'''+15y^2y''^3+9yy'^2y''^2+3y'^4y''=0~,\label{example2}
\end{eqnarray} is of the form (\ref{fourthindx-lform-type1}) with the coefficients
$A_1={-10\over y'}$, $A_0={-3y'\over y}$, $B_3={15\over {y'^2}},$ $B_2={9\over y},$ $B_1={3y'^2\over
{y^2}},$ $B_0=0$
satisfying the conditions (\ref{fourthindx-lc-type2}). So it is reduced to the third order
linearizable ODE
\begin{eqnarray}
y^2u^2u'''-3yu^2u''-6y^2uu'u''+3u^2u'+6yuu'^2+6y^2u'^3=0~,\label{reduced8}
\end{eqnarray}
with $y$ as independent and $u$ as dependent variable. The transformation $t=y^2,$ $s={1\over
u},$ reduces (\ref{reduced8}) to the linear third order ODE
$s'''=0,$ whose solution is $s=c_1t^2+c_2t+c_3.$ Now one only needs to solve the
equation $y'=1/(c_1y^4+c_2y^2+c_3)$, where $c_i$ are arbitrary constants. Hence, the general solution of (\ref{example2}) is given by
\begin{eqnarray}
x=c_1y^5+c_2y^3+c_3y+c_4~.\nonumber
\end{eqnarray}
\\
{\bf Example 3.} The ODE
\begin{eqnarray}
y'y''y^{(4)}-3y'y'''^2+6y'^3y''^2y'''-4y''^2y'''-y'y''^5=0~,\label{example3}
\end{eqnarray}
has 2 symmetries. It is of the form (\ref{fourthindx-lform-type2}) with the coefficients $r_0=0,$ $C_2=6y'^2-{4\over y'},$ $C_1=C_0=0,$ $D_5=-1,$ $D_4=D_3=D_2=D_1=D_0=0,$
obey the conditions (\ref{C1})--(\ref{C9}). So it is reducible to linearizable
third order ODE
\begin{eqnarray}
u'''+{1\over u'}[-3u''^2-yu'^5]=0~.\label{reduced9}
\end{eqnarray}
The transformation
$t=u,s=y,$
will convert the nonlinear ODE (\ref{reduced9}) to the linear ODE
$s'''+s=0$ with solution
\begin{eqnarray}
s=c_1e^{-t}+c_2e^{t\over 2}\cos{t}+c_3e^{t\over 2}\sin{t}~.
\end{eqnarray}
Finally to find the solution of (\ref{example3}), we only need to solve
\begin{eqnarray}
y=c_1e^{-y'}+c_2e^{y'\over 2}\cos{y'}+c_3e^{y'\over 2}\sin{y'}~.
\end{eqnarray}
\\
\textbf{Example 4.} The nonlinear ODE
\begin{eqnarray}
y''y^{(4)}+y'''^3-y'''^2-y''y'''~,\label{example4}
\end{eqnarray}
is of the form (\ref{fourthindxy-lform}) and the coefficients
$a={1\over y''},$ $b=-{1\over y''},$ $c=-{y''},$ $d=0$,
that satisfy conditions (\ref{fourthindxy-lc}). So it is reduced to the linearizable
second order ODE $u''+u'^3-u'=0$.
By using the transformation $t=u~,~s=e^y~$, we can reduce it to
linear ODE
$s''-s=0$,
whose solution is given by
$s=c_1e^t+c_2e^{-t}~,$
where $c_i$ are arbitrary constants. So that solution of (\ref{example4}) is
obtained by solving the second order ODE
\begin{eqnarray}
e^{y'}=c_1e^{-y''}+c_2e^{y''}~,
\end{eqnarray}
where $c_i$ are arbitrary constants.

\section{Concluding Remarks}
Nonlinear ODEs are difficult to solve but, if they can be
converted to linear ones by invertible transformations, they can
be solved. Hence linearization plays a significant role in the
theory of ODEs. In this paper we have presented criteria for fourth
order autonomous ODEs to be reducible to linearizable third and
second order ODEs. There are certain fourth order ODEs, not
depending explicitly on the independent variable, which cannot be
linearized by point or contact transformations but can be reducible
to linearizable third order ODEs by Meleshko's method. The
solution of the original equation is then obtained by a quadrature.
Various fourth order ODEs with fewer symmetries can be reduced to
linearizable form by this procedure. The class of ODEs linearizable
by this method is not included in the Ibragimov and Meleshko
classes or conditionally linearizable classes \cite{mq3, mq4} of the ODEs (though
there can be an overlap but it is not contained in that either). The
reason is that it is not linearizable but reducible to linearizable
form. In Lie's programme there is no definite statement available
for the cases when the ODEs are not linearizable. By the recent
developments this gap may be filled. By using the concept of
Meleshko linearization a new class of scalar ODEs may be defined on
the basis of initial
conditions to be satisfied by ODEs.\\
\begin{center}
\textbf{Appendix}
\end{center}
\begin{eqnarray}
(r_0C_1-6r_{0y})y'^2+(6r_0r_{0y'}+4r_0^2-r_0^2C_2-C_0)y'-4r_0^2=0~,\label{C1}\\
\nonumber\\
(C_{2y}-C_{1y'})y'^3+(r_0C_{2y'}+C_2r_{0y'}-4r_{0y'}-6r_{0y'y'})y'^2\nonumber\\
+(10r_{0y'}+4r_0-C_2r_0)y'-8r_0=0~,\label{C2}\\
\nonumber\\
(-6r_0^2C_{1y}-54{(r_{0y})}^2+18r_0r_{oyy}+18r_0r_{0y}C_1-2r_0^2C_1^2)y'^8\nonumber\\
+(3r_0^3C_{1y'}+48r_0^2r_{0y}-3r_0^3C_{2y}-36r_0^2r_{0yy'}
-6r_0^2r_{0y}C_2-18r_0^2r_{0y'}C_1\nonumber\\
+2r_0^3C_1C_2-16r_0^3C_1)y'^7+(-60r_0^3r_{0y'}+9r_0^4C_{2y'}-42r_0^2r_{0y}\nonumber\\
-36r_0^2{(r_{oy'})}^2+9r_0^3r_{0y'}C_2+14r_0^3C_1-32r_0^4+8r_0^4C_2+4r_0^4C_2^2\nonumber\\
+18r_0^4D_4)y'^6+(44r_0^4+72r_0^2r_{0y'} -18r_0^3r_{0y'}
-7r_0^4C_2)y'^5\nonumber\\
+(-20r_0^4)y'^4-72r_0^5D_5=0~,\label{C3}\\
\nonumber\\
(-12r_0C_{1y}+18r_{oyy'}+18r_{0y}C_1-4r_0C_1^2)y'^8+(9r_0^2C_{1y'}-48r_0r_{0y}\nonumber\\
-27r_0^2C_{2y}-36r_0r_{0yy'}-18r_{0y}+72r_0r_{0y}+24r_0r_{0y}C_2-18r_0r_{0y'}C_1\nonumber\\
-18r_0r_{0y'}-32r_0^2C_1-2r_0^2C_1C_2)y'^7+(-18D_1-36r_0^2r_{0y'}+33r_0^3C_{2y'}\nonumber\\
+6r_0r_{0y}+18r_0^2C_1-21r_0^2r_{0y'}C_2+18r_0{(r_{0y'})}^2-64r_0^3+4r_0^2C_1-8r_0^3C_2\nonumber\\
+20r_0^3C_2^2+72r_0^3D_4)y'^6+(52r_0^3+6r_0^2r_{oy'}+13r_0^3c_2)y'^5\nonumber\\
+(-22r_0^3)y'^4-270r_0^4D_5=0~,\label{C4}\\
\nonumber\\
(-3C_{1y}-C_1^2)y'^8+(3r_0C_{1y'}-12r_{0y}-21r_0C_{2y}-8r_0C_1\nonumber\\
+15r_{oy}C_2-5r_0C_1C_2)y'^7+(-9d_2+12r_0r_{0y'}+21r_0^2C_{2y'}-30r_{0y}\nonumber\\
-15r_0r_{0y'}C_2+10r_0C_1-20r_0^2C_2+14r_0^2C_2^2+54r_0^2D_4\nonumber\\
-16r_0^2)y'^6+(-9C_0+28r_0^2+30r_0r_{0y'}+13r_0^2C_2)y'^5+(-40r_0^2)y'^4\nonumber\\
-180r_0^3D_5=0~,\label{C5}\\
\nonumber\\
(-3C_{2y}-C_1C_2)y'^7+(-3D_3+4C_1+3r_0C_{2y'}-4r_0C_2+2r_0C_2^2\nonumber\\
+12r_0D_4)y'^6+(-4r_0+4r_0C_2)y'^5+(-r_0)y'^4-30r_0^2D_5=0~,\label{C6}\\
\nonumber\\
(-54D_{4y}+18C_{1y'y'}+3C_2C_{1y'}-72C_{2yy'}-39C_2C_{2y})y'^8+(24C_{2y}\nonumber\\
+72r_{0y'y'}+12C_2r_{0y'}-6C_{1y'}+36r_0C_{2y'y'}-3r_0C_2C_{2y'}+72r_{0y'}C_{2y'}\nonumber\\
+33C_2^2r_{0y')}+108D_4r_{0y'}+54r_0d_{4y'}+36r_0C_2^2+18r_0C_{2y'y'})y'^7\nonumber\\
+(-168r_{0y'}-12r_0C_2-138r_0C_{2y'}-24C_2r_{0y'}-33r_0C_2^2-36r_0D_4)y'^6\nonumber\\
+(168r_0-228r_0C_2+60r_{0y'})y'^5+(-120r_0)y'^4+(270D_5r_{0y}\nonumber\\
+270r_0D_{5y})y'^2+(54r_0^2D_{5y'}-810r_0r_{0y'}D_5)y'+2160r_0^2D_5=0~,\label{C7}
\end{eqnarray}
and
\begin{eqnarray}
(-H_y)y'^2+(3Hr_{0y'}+r_0H_y')y'-3Hr_0=0~,\label{C8}
\end{eqnarray}
where
\begin{eqnarray}
H=(D_{4y'}+{1\over 3}C_{2y'y'}+{2\over 3}C_2C_{2y'}+{2\over
3}C_2D_4+{4\over 27}C_2^3)\nonumber\\
+{1\over y'}(-{4\over 3}C_{2y'}+{2\over 3}C_2^2-{4\over
3}D_4-{8\over 9}C_2^2)\nonumber\\
+{1\over y'^2}(-{5\over
9}C_2)+{1\over y'^3}({40\over 27})+{1\over y'^5}(-2D_{5y}-{2\over
3}C_1D_5)\nonumber\\
+{1\over y'^6}(-3r_0D_{5y'}-5D_5r_{0y'}-2r_0C_2D_5-{8\over
3}r_0D_5)\nonumber\\
+{1\over y'^7}(24r_0D_5)~.\label{C9}
\end{eqnarray}


\begin{thebibliography}{99}
\bibitem{fm1}
Mahomed FM, Point symmetry group classification of ordinary
differential equation: A survey of some results,
\textit{Mathematical Methods in the Applied Sciences} {\bf 30}
(2007) 1995-2012.
\bibitem{l1}
Lie S, Theorie der transformationsgruppen, \textit{Math. Ann.} {\bf
16} (1880) 441.
\bibitem{l2}
Lie S, Klassifikation und Integration von gew\"{o}nlichen
Differentialgleichungenzwischen $x$, $y$, die eine Gruppe von
Transformationen gestaten, \textit{Arch. Math} {\bf VIII, IX} (1883)
187.
\bibitem{t} Tresse A. Sur les invariants differentiels des
groupes continus de transformations. \textit{Acta Math}. {\bf 18}
(1894) 1.
\bibitem{c1} Chern SS. Sur la g\'{e}om\'{e}trie d'une \'{e}quation
diff\'{e}rentielle du trois\`{e}me ordre. \textit{CR Acad Sci Paris}
(1937) 1227-1229.
\bibitem{c2} Chern SS. The geometry of the differential
equation $y''''=F(x,y,y'',y''')$. \textit{Sci Rep Nat Tsing Hua
Univ} {\bf. 4} (1940) 97-111.

\bibitem{g1} Grebot G. The linearization of third order
ODEs, preprint (1996).
\bibitem{g2} Grebot G. The characterization of third order
ordinary differential equations admitting a transitive
fibre-preserving point symmetry group. \textit{J Math Anal Appl}
{\bf 206} (1997) 364-388.
\bibitem{mL} Mahomed FM, Leach PGL. Symmetry Lie algebras of
nth order ordinary differential equations. \textit{J Math Anal Appl}
{\bf 151} (1990) 80-107.
\bibitem{np}
Neut S, Petitot M. La g\'{e}om\'{e}trie de l'\'{e}quation
$y'''=f(x,y,y',y'')$. \emph{CR Acad Sci Paris S\'{e}r I}
\textbf{335} (2002) 515-518.
\bibitem{im}
Ibragimov NH, Meleshko SV. Linearization of third order ordinary
differential equations by point and contact transformations.
\textit{J Math Anal Appl} {\bf 308} (2005) 266-289.
\bibitem{ims1}
Ibragimov NH, Meleshko SV, Suksern S. Linearization of fourth order
ordinary differential equation by point transformations. \textit{J
Phys. A: Math. Theor.} {\bf 41} (2008).
\bibitem{ims2}
Suksern S, Ibragimov NH, Meleshko SV. Criteria for the fourth order
ordinary differential equations to be linearizable by contact
transformations. \textit{Common Nonlinear Sci Number Simulat} {\bf
14} (2009) 2619-2618.
\bibitem{M} Meleshko SV. On linearization of third order
ordinary differential equation. \textit{J Phys A Math Gen} {\bf 39}
(2006) 15135-15145.
\bibitem{mq3} Mahomed FM, Qadir A. Conditional linearizability criteria for third order ordinary differential equations. \emph{J. Nonlinear Math. Physics.} \textbf{15}(Suppl. 1) (2008) 124-133.
    \bibitem{mq4} Mahomed FM, Qadir A. Conditional linearizability of fourth order semi-linear ordinary differential equations. \emph{J. Nonlinear Math. Physics.} \textbf{16} (2009) 165-178.
        \bibitem{mq2} Mahomed FM, Qadir A. Classification of ordinary
differential equations by using conditional linearizability and
symmetry. \emph{Commun. Nonlinear Sci. Numer. Simulat.} \textbf{17}
(2012) 573-584.


\end{thebibliography}
\end{document}